\documentclass{article}
\usepackage{amsmath,graphicx,amsthm,geometry}
\usepackage[all]{xy}

\newtheorem{conjecture}{Conjecture}
\newtheorem{lemma}{Lemma}
\newtheorem{proposition}{Proposition}

\newtheorem{remark}{Remark}

\newcommand{\defn}{\ensuremath{\overset{\mathrm{def}}{=}}}

\newcommand{\dd}{\ensuremath{\mathrm{d}}}

\begin{document}
\title{\bf A Nonlocal Formulation of\\ Rotational Water Waves}
\author{A. C. L. Ashton and A. S. Fokas}
\date{4th March 2011}
\maketitle

\begin{abstract}
The classical equations of irrotational water waves have recently been reformulated as a system of two equations, one of which is an explicit non-local equation for the wave height and for the velocity potential evaluated on the free surface. Here, in the two dimensional case, (a) we generalise the relevant formulation to the case of constant vorticity, as well as to the case where the free surface is described by a multi-valued function; (b) in the case of travelling waves we derive an upper bound for the free surface; (c) in the case of constant vorticity we construct a sequence of nearly Hamiltonian systems which provide an approximation in the asymptotic limit of certain physical small parameters. In particular, the explicit dependence of the vorticity on the coefficients of the KdV equation is clarified.  Also, in the irrotational case we extend the formalism to $n>2$ dimensions and analyse rigorously the linear limit of these equations.
\end{abstract}

\section{Introduction}
We consider the classical problem in hydrodynamics concerning the propagation of surface waves generated by an incompressible fluid with free surface. In the case of irrotational flow we consider the problem in $n=2,3$ (or higher) spatial dimensions, whereas in the rotational case we confine attention to $n=2$. We consider the problem in its full generality so that the depth of the fluid may not be constant, and consider \emph{all} solutions, not only those that arise under the travelling wave assumption. In addition, we include the effect of surface tension so that capillary-type waves are included in our study.

We denote by $\mathcal{B}_h$ the bottom surface:
\begin{equation} \mathcal{B}_h = \{ (x,y): x\in \mathbf{R}^{n-1},\, y= -h_0 + h(x) \}, \end{equation}
where $h_0$ is constant and $h(x)$ is a real valued function. We denote by $\mathcal{S}_\eta$ the free surface:
\begin{equation} \mathcal{S}_\eta  =\{ (x,y): x\in\mathbf{R}^{n-1},\, y=\eta(x,t)\} \qquad \textrm{for $t\geq 0$.} \end{equation}
We refer to $\eta$ as the height of the wave and we assume $\eta+h_0>h$ for each $x\in\mathbf{R}^{n-1}$. The domain of the problem is the region between $\mathcal{B}_h$ and $\mathcal{S}_\eta $, which is denoted by $\Omega$.

\subsection{The Irrotational Case in $n\geq 2$ Dimensions}
In the irrotational case we introduce the velocity potential $\phi$, where $u=\nabla\phi$, and then the governing equations become:
\begin{subequations}\label{irrotational-n}
\begin{alignat}{2}
\Delta \phi &= 0 &\qquad  &\textrm{in $\Omega$,} \label{e1}\\
\nabla\phi \cdot N_\mathcal{B} &=0 &&\textrm{on $\mathcal{B}_h$,}\label{e2} \\
\nabla\phi \cdot N_\mathcal{S} &=\eta_t  && \textrm{on $\mathcal{S}_\eta $,} \label{e3} \\
\phi_t + \tfrac{1}{2} |\nabla\phi |^2 + g\eta &= f(\eta) && \textrm{on $\mathcal{S}_\eta $,}\label{e4}
\end{alignat}
\end{subequations}
where $g$ is the acceleration due to gravity, $\nabla$ denotes the usual vector gradient, $N_\mathcal{B}$ is the exterior normal to $\mathcal{B}_h$ and $N_\mathcal{S}$ is the exterior normal to $\mathcal{S}_\eta$, i.e.
\begin{subequations}\label{defs}
\begin{equation}\nabla = (\nabla \! _x, \partial_y),\qquad N_\mathcal{B} =(-\nabla \! _x h, -1),\qquad N_\mathcal{S}=(-\nabla \! _x\eta,1), \label{normals}\end{equation}
where $\nabla \! _x=(\partial_{x_1}\ldots,\partial_{x_{n-1}})$ is the horizontal gradient. The right hand side of \eqref{e4} is functionally dependent on $\eta$ through:
\begin{equation} f(\eta) = \frac{\sigma}{\rho} \nabla \! _x \cdot\left( \frac{ \nabla \! _x \eta}{\sqrt{1+|\nabla \! _x \eta|^2}}\right), \label{stens} \end{equation}
\end{subequations}
which is a measure of the effect of surface tension ($\sigma$ and $\rho$ denote the constant surface tension and density respectively). Equation \eqref{e1} is a consequence of incompressibility, \eqref{e2} is the Neumann condition satisfied on $\mathcal{B}_h$ and \eqref{e3}, \eqref{e4} are the kinematic, Bernoulli conditions respectively on $\mathcal{S}_\eta$.

We will reformulate the problem in terms of the functions $\eta$ and $q$ where $q(x,t) = \phi(x,\eta,t)$, i.e. $q$ represents the potential on the free surface $\mathcal{S}_\eta $. It was shown by Zakharov in his classical paper (\cite{zakharov1968spw}) that in the case $n=2$ the pair $(\eta,q)$ constitutes a canonically conjugate pair in the Hamiltonian formulation of this problem. The generalisation of Zakarov's result to higher dimensions is straightforward. Indeed, in the case of constant $\mathcal{B}_h$ and zero surface tension equations \eqref{irrotational-n} admit a Hamiltonian formulation with respect to the Hamilontian
\begin{equation} H= \iint_\Omega \tfrac{1}{2} |\nabla\phi |^2\, \dd x\, \dd y + \int \tfrac{1}{2} g\eta^2\, \dd x\label{hamiltonian} \end{equation}
and the standard symplectic structure. Here and throughout $\dd x = \dd x_1 \cdots  \dd x_{n-1}$ denotes the Lebesgue measure on $\mathbf{R}^{n-1}$.

Applying the chain rule to the expression $q(x,t)=\phi (x,\eta,t)$ we find the following relations that are valid on free surface $y=\eta(x,t)$:
\begin{subequations}\label{chain}
\begin{align}
\nabla \! _x q &= \nabla \! _x \phi + (\partial_y \phi) \nabla \! _x \eta ,\label{chainsp} \\
q_t &= \phi_t + \phi_t \eta_t. 
\end{align}
\end{subequations}
It is possible to solve for $\nabla\phi$ in terms of the $q$ and $\eta$ following a similar calculation with that of \cite{ablowitz2006nnl}. Indeed, using \eqref{e3} in equations \eqref{chainsp} we find the following nonsingular set of equations for $\nabla \! _x \phi$:
\begin{equation}\nabla \! _x q - \eta_t \nabla \! _x \eta = (\mathbf{I} + \nabla \! _x \eta \otimes \nabla \! _x \eta) \cdot \nabla \! _x \phi \end{equation}
and then equation \eqref{e3} gives $\phi_y$ in terms of $q$ and $\eta$. In the new variables the dynamic boundary condition on $\mathcal{S}_\eta $, i.e. \eqref{e4}, becomes
\begin{equation}
q_t + \tfrac{1}{2} |\nabla \! _x q|^2 + g\eta - \frac{(\eta_t + \nabla \! _x \eta \cdot \nabla \! _x q)^2}{2(1+|\nabla \! _x \eta |^2)} = f(\eta). \label{e5}
\end{equation}
In what follows, the notation $\int \dd x$ will be used to denote an integral over $\mathbf{R}^{n-1}$.

\subsubsection{Bounds on Wave Height}
Let us temporarily confine attention to travelling wave solutions to \eqref{bernoullirot}, so that $\eta=\eta(x-ct)$ and $\xi=\xi(x-ct)$. Setting $z=x-ct$ and denoting differentiation with respect to $z$ with a prime, the equation \eqref{bernoullirot} becomes
\begin{equation}
-c \xi' + \tfrac{1}{2} (\xi')^2 + (g+\gamma c)\eta + \frac{\gamma\eta (-2c(\eta')^2 - 2\xi'+\gamma\eta)}{2(1+(\eta')^2)} - \frac{(-c\eta' + \eta' \xi')^2}{2(1+(\eta')^2)}  = f(\eta). \label{g2''}
\end{equation}
Following \cite{tio2010}, we observe \eqref{g2''} is actually a quadratic in the unknown $\xi '$, which can be solved to give
\begin{equation} \xi'(z) = \gamma \eta + c - \sqrt{ \left(1+(\eta')^2\right)\left( c^2 + f(\eta) -2g\eta \right)}, \label{solve_quad'}\end{equation}
where we have discarded the second solution which fails to decay for large $z$. 

Equation \eqref{solve_quad'} immediately gives the estimate:
\begin{equation} \eta \leq  \frac{c^2}{2g} + \frac{\sigma}{2\rho g} \frac{\dd}{\dd z} \left( \frac{\eta '}{\sqrt{1+(\eta')^2}} \right). \label{upperboundheight} \end{equation}
In the absence of vorticity and surface tension, this reproduces the upper bound derived in \cite{craig1988symmetry}. With the assumption that $|\eta'|$ is bounded, we can integrate this inequality to establish the following estimate for the average wave height:
\[ \langle \eta\rangle \leq \frac{c^2}{2 g},  \]
which holds in the presence of arbitrary surface tension.

\subsection{The Rotational Case in Two Dimensions}\label{rot}
In this section we derive the equations for the free surface problem for a two dimensional fluid with constant vorticity. Denoting the velocity of the flow by $(u,v)$, the Euler equations for invisid flow are the following equations:
\begin{subequations}
\begin{align}
u_t + uu_x + vu_y &= -P_x, \label{Euler1} \\
v_t + uv_x + vv_y &= -P_y - g, \label{Euler2} \\
u_x + v_y &=0, \label{incomp}
\end{align}
\end{subequations}
where $P=P(x,y)$ is the fluid pressure. Let $\omega$ denote the vorticity of the fluid, i.e.
\begin{equation} \omega = v_x - u_y. \label{vorticity}\end{equation}
Eliminating $P$ from \eqref{Euler1} and \eqref{Euler2} we find:
\begin{equation}
\frac{\partial \omega}{\partial t} + (u\partial_x + v\partial_y)\omega = 0. \label{vort}
\end{equation}
We restrict attention to the case in which $\omega=\gamma$ is constant throughout $\Omega$, so that \eqref{vort} is satisfied identically. The domain $\Omega$ is simply connected, therefore we can introduce a globally defined stream function $\psi(t,x,y)$ so that:
\begin{equation} u= \psi_y, \qquad v = -\psi_x, \qquad (x,y)\in\Omega. \end{equation}

Replacing in equation \eqref{vorticity} $u$ by $\psi_y$ and $v$ by $-\psi_x$ we find $\Delta \psi = -\gamma$. Thus, the function defined by $\psi^h = \psi + \tfrac{1}{2} y^2\gamma$, is harmonic in $\Omega$. In fact, $\psi^h$ is only defined up to the addition of an arbitrary function of time, so by abuse of notation we write $\psi^h$ to mean $[\psi^h]$, i.e. the equivalence class of such functions that differ only by a function of time. Let $\varphi$ be the harmonic conjugate of $\psi^h$, i.e.
\begin{subequations}
\begin{alignat}{2}
\varphi_x &= \psi_y + \omega y, \label{varphi_x} \\
\varphi_y &= -\psi_x. \label{varphi_y}
\end{alignat}
\end{subequations}
The function $\varphi$ is harmonic throughout $\Omega$.

Adding and subtracting in equations \eqref{Euler1} and \eqref{Euler2} the terms $vv_x$ and $uu_y$ respectively, the equations become
\begin{subequations}\label{Eulernew}
\begin{align}
u_t + \tfrac{1}{2} \partial_x (u^2+v^2) -\gamma v   &= -P_x, \\
v_t + \tfrac{1}{2} \partial_y (u^2+v^2) + \gamma u &= -P_y-g. 
\end{align}
\end{subequations}
Using in equations \eqref{Eulernew} the identities
\[ u_t = (\psi_y)_t = (\varphi_x)_t, \quad v= -\psi_x, \quad v_t = -(\psi_x)_t = (\varphi_y)_t, \quad u=\psi_y ,\]
equations \eqref{Eulernew} can be integrated to give the following equation:
\begin{equation} \varphi_t + \tfrac{1}{2} |\nabla\psi|^2 + \gamma \psi + P + gy = \alpha(t), \qquad (x,y)\in \Omega, \end{equation}
where $\alpha(t)$ is some function of time. Since we are still dealing with an equivalence class of functions $\psi$ and $\varphi$ we can absorb a function of time into them, so we arrive at
\begin{equation} \varphi_t + \tfrac{1}{2} |\nabla\psi|^2 + \gamma \psi + gy = P_{\mathrm{atm}}-P,\qquad (x,y)\in\Omega, \label{bc}\end{equation}
where $P_\mathrm{atm}$ is the atmospheric pressure above the free surface $\mathcal{S}_\eta$. On the free surface $\mathcal{S}_\eta$ we have the dynamic boundary condition
\begin{equation} P_\mathrm{atm} -P = f(\eta), \qquad f(\eta) \defn \sigma  \left(\frac{\eta_x}{\sqrt{1+\eta_x^2}}\right)_x. \label{bcc}\end{equation}
The stream function $\psi$ is \emph{not} constant on the free surface. In particular if we introduce the function $\Psi(x,t)=\psi(x,\eta(x,t),t)$ to be the stream function defined on the free surface, we see
\begin{align*}
\Psi_x &= \psi_x + \eta_x \psi_y \\
 &= -v + u\eta_x \\
 &= -\eta_t,
\end{align*}
which follows from the kinematic boundary condition \eqref{e3} on the free surface. So up to an additional function of time, which again we simply absorb into the definition of $\psi$, we have
\[ \psi|_\mathcal{S} = \Psi(x,t) = -\int_{-\infty}^x \eta_t(y,t)\, \dd y \equiv -\partial_x^{-1} \eta_t. \]
Hence \eqref{bc} and \eqref{bcc} imply the following nonlinear boundary condition:
\begin{equation}
\partial_t \varphi- \gamma\partial_x^{-1} \eta_t + \tfrac{1}{2} |\nabla\psi|^2 + g\eta =f(\eta) \qquad \textrm{on $\mathcal{S}_\eta$}. \label{newbc}
\end{equation}
The kinematic condition on $\mathcal{S}_\eta$ and the condition on $\mathcal{B}_h$ are the same as in the irrotational case. 

\emph{In summary, if the vorticity equals the constant $\gamma$, then the functions $\varphi(x,y,t)$ and $\eta(x,t)$ satisfy the following boundary value problem:}
\begin{subequations}\label{roteqns}
\begin{alignat}{2}
\Delta \varphi &= 0 &\qquad  &\textrm{in $\Omega$,} \label{f1}\\
(\varphi_x - \gamma y, \varphi_y ) \cdot N_\mathcal{B} &=0 &&\textrm{on $\mathcal{B}_h$,}\label{f2} \\
(\varphi_x - \gamma y, \varphi_y) \cdot N_\mathcal{S} &= \eta_t && \textrm{on $\mathcal{S}_\eta $,} \label{f3} \\
\varphi_t - \gamma\partial_x^{-1} \eta_t+ \tfrac{1}{2} \varphi_y^2 + \tfrac{1}{2}(\varphi_x - \gamma y)^2 + gy &= f(\eta) && \textrm{on $\mathcal{S}_\eta $,}\label{f4}
\end{alignat}
\end{subequations}
where $N_\mathcal{B}$, $N_\mathcal{S}$ and $f(\eta)$ are defined in \eqref{defs}.

As in the irrotational case, we will formulate the problem in terms of $\eta$ and $\xi$ where $\xi(x,t)=\varphi(x,\eta,t)$, i.e $\xi$ represents the pseudo-potential $\varphi$ on the free surface. Computing $\xi_x$ using the chain rule (compare with \eqref{chainsp}) and invoking \eqref{f3} we find the following pair of equations for $(\varphi_x, \varphi_y)$ on $\mathcal{S}_\eta$:
\begin{align*}
\xi_x &= \varphi_x + \eta_x \varphi_y, \\
\eta_t - \gamma \eta \eta_x &= \varphi_y - \eta_x \varphi_x.
\end{align*}
Solving these equations for $(\varphi_x, \varphi_y)$ we find the following relations on $\mathcal{S}_\eta$:
\begin{subequations}\label{mess}
\begin{align}
(1+\eta_x^2) \varphi_x &= \xi_x - \eta_t \eta_x + \gamma \eta \eta_x^2, \\
(1+\eta_x^2) \varphi_y &= \eta_t + \eta_x \xi_x - \gamma\eta \eta_x.
\end{align}
\end{subequations}
Furthermore, using \eqref{mess} in the boundary condition \eqref{f4} on $\mathcal{S}_\eta$ we find:
\begin{equation}
\xi_t  + \tfrac{1}{2} \xi_x^2 + g\eta - \gamma\partial_x^{-1} \eta_t + \frac{\gamma\eta (2\eta_t \eta_x - 2\xi_x+\gamma \eta)}{2(1+\eta_x^2)} - \frac{(\eta_t + \eta_x \xi_x)^2}{2(1+\eta_x^2)}  = f(\eta). \label{bernoullirot}
\end{equation}
It is reassuring to note that \eqref{bernoullirot} reduces to \eqref{e5} when $\gamma=0$. 

\subsection{Multivalued Free surface}
In the above analysis we have assumed that $\eta(x,t)$ is a single valued function. However, in general, this assumption can be violated. In particular, it is known that in the case of water waves with constant vorticity, $\eta(x,t)$ quickly becomes a multi-valued function, as observed by \cite{telesdasilva1988sss}. In what follows, we show that it is conceptually straightforward to modify our analysis so that it can address the case of a multi-valued free surface.

The free surface is assumed to be a one-dimensional $C^2$-differentiable manifold, so features such as cusps and self-intersections are ignored. We set:
\[ \mathcal{S}= \{ (x,y)\in \mathbf{R}^2: x=X(\lambda,t), y=Y(\lambda,t) \} ,\qquad t\geq 0 \]
where $\lambda\in\mathbf{R}$ is the arc-length along the curve and $X(\lambda, \cdot)$, $Y(\lambda, \cdot)$ are assumed to be $C^2$. The outward normal is now given by $N_\mathcal{S}=(-\dot{Y},\dot{X})$, where the dot denotes differentiation with respect to $\lambda$. The kinematic boundary condition in \eqref{f3} now becomes:
\begin{equation}(\varphi_x - \gamma y - X_t , \varphi_y-Y_t) \cdot N_\mathcal{S} =0, \qquad \textrm{on $\mathcal{S}$.} \label{newkin}\end{equation}
Again introducing $\Psi(\lambda,t)=\psi(X(\lambda,t),Y(\lambda,t),t)$ to be the stream function evaluated on the free surface, the chain rule and the new kinematic boundary condition \eqref{newkin} gives
\[ \dot{\Psi} = \dot{X} Y_t - \dot{Y}X_t \equiv \frac{\partial(X,Y)}{\partial (\lambda,t)}. \]
Using analogous reasoning as in the case of the single valued free surface, this gives
\begin{equation} \psi|_\mathcal{S} = - \partial_\lambda^{-1} \left( \dot{X} Y_t - \dot{Y}X_t \right). \label{multtemp}\end{equation}
Again, we introduce a function $\xi$ to represent the pseudo-potential on the free surface, i.e. $\xi(\lambda,t) = \varphi(X(\lambda,t),Y(\lambda,t),t)$. In analogy with equations \eqref{chain}, we note the following relations on $\mathcal{S}$:
\begin{subequations}\label{newchain}
\begin{align}
 \dot{\xi} &= \dot{X} \varphi_x + \dot{Y} \varphi_y,\\
 \xi_t &= X_t \varphi_x + Y_t \varphi_y + \varphi_t.
\end{align}
\end{subequations}
Equations \eqref{newkin} and \eqref{newchain} constitute a non-singular system for the functions $\{\varphi_x, \varphi_y,\varphi_t\}$ on the free surface $\mathcal{S}$, so we can invert these relationships to express these functions in terms of $\{X,Y, \eta\}$. After some algebra, one finds the following relations on $\mathcal{S}$:
\begin{subequations}\label{inverted}
\begin{align}
\varphi_x &= \frac{\dot{\xi} \dot{X} + \dot{Y} (\dot{Y}(\gamma Y + X_t)-\dot{X}Y_t)}{(\dot{X}^2+\dot{Y}^2) } \\
\varphi_y &= \frac{\dot{\xi} \dot{Y} - \dot{X} (\dot{Y}(\gamma Y+X_t)-\dot{X}Y_t)}{(\dot{X}^2+\dot{Y}^2) } \\
\varphi_t &= \xi_t  - \frac{\dot{\xi}(\dot{X}X_t + \dot{Y}Y_t)+ (X_t \dot{Y} - Y_t\dot{X})(\dot{Y}(\gamma Y + X_t)-\dot{X}Y_t)}{(\dot{X}^2+\dot{Y}^2) }.
\end{align}
\end{subequations}
Using these expressions on $\mathcal{S}$ and the relation in \eqref{multtemp}, the left hand side of the Bernoulli condition in \eqref{f4} becomes:
\begin{multline*} \xi_t + gY -\tfrac{1}{2}X_t^2 - \tfrac{1}{2}Y_t^2 - \gamma \partial_\lambda^{-1} \left( \dot{X} Y_t - \dot{Y}X_t \right) \\
+ \frac{\gamma Y(2\dot{X}\dot{Y}Y_t -2\dot{X}^2X_t-2\dot{\xi}\dot{X}+\gamma \dot{X}^2Y)}{2(\dot{X}^2+\dot{Y}^2)} + \frac{(\dot{\xi}-\dot{X}X_t-\dot{Y}Y_t)^2}{2(\dot{X}^2+\dot{Y}^2)}
\end{multline*}
and the term due to surface tension is:
\[ \frac{\sigma}{\rho} \frac{\dot{X}\ddot{Y}-\dot{Y}\ddot{X}}{(\dot{X}^2+\dot{Y}^2)^{3/2}}. \]
Here the coefficient of $\sigma$ is the intrinsic curvature of the surface $\mathcal{S}$. Combining these expressions gives the Bernoulli condition on $\mathcal{S}$ in terms of the surface parameters $(X,Y)$ and the potential on the free surface, $\xi(\lambda,t)$:
\begin{multline} \xi_t + gY -\tfrac{1}{2}X_t^2 - \tfrac{1}{2}Y_t^2 - \gamma \partial_\lambda^{-1} \left( \dot{X} Y_t - \dot{Y}X_t \right)\\ + \frac{\gamma Y(2\dot{X}\dot{Y}Y_t -2\dot{X}^2X_t-2\dot{\xi}\dot{X}+\gamma \dot{X}^2Y)}{2(\dot{X}^2+\dot{Y}^2)} + \frac{(\dot{\xi}-\dot{X}X_t-\dot{Y}Y_t)^2}{2(\dot{X}^2+\dot{Y}^2)} - \sigma \frac{\dot{X}\ddot{Y}-\dot{Y}\ddot{X}}{(\dot{X}^2+\dot{Y}^2)^{3/2}} =0. \label{bigbernoulli}
\end{multline} 
We note that in the case $X=\lambda$, $Y=Y(X,t)$ this equation reduces to \eqref{bernoullirot}.

\subsubsection{Bounds on Wave Height}
One can also consider travelling wave solutions in this setting and achieve an analogous upper bound to \eqref{upperboundheight}. In this case
\begin{align*}
 X(\lambda,t)&\mapsto X(\lambda) -c t \\
 Y(\lambda,t)&\mapsto Y(\lambda).
\end{align*}
Setting $\xi = \varphi (X-ct,Y)$ with $X=X(\lambda)$ and $Y=Y(\lambda)$ being independent of time, equation \eqref{bigbernoulli} reduces to a quadratic in $\dot{\xi}$. Examining the discriminant of this quadratic reveals the estimate
\[ Y \leq \frac{c^2}{2g} + \frac{\sigma}{2\rho g}\left[ \frac{\dot{X}\ddot{Y}-\dot{Y}\ddot{X}}{(\dot{X}^2+\dot{Y}^2)^{3/2}}\right]. \]
We deduce that the estimate in \eqref{upperboundheight} for the wave height is valid for \emph{all} travelling waves, including those which are truly multi-valued.

The above analysis shows that although it is straightforward to incorperate the case of a multivalued free surface, the relevant formulae become more complicated, so for convenience we will present most of our results assuming that $\eta(x,t)$ is single valued.

\section{The Non-Local Formulation}
A novel non-local formulation in two and three dimensions for irrotational water waves was presented by \cite{ablowitz2006nnl}. Here we extend this formulation to rotational water waves with constant vorticity in two dimensions and for irrotational water waves in an arbitrary number of dimensions.

The formulation by \cite{ablowitz2006nnl} is based on the general approach for studying boundary value problems for a large class of partial differential equations introduced by one of the authors, see \cite{fokas1997utm} and \cite{fokas2000integrability}. A crucial role in this approach is played by a certain nonlocal equation called the \emph{global relation}. The global relation is a consequence of the following fact: suppose that the functions $u$ and $v$ are harmonic in $\Omega\subset \mathbf{R}^{n-1}\times\mathbf{R}$. Then,
\begin{equation}
\nabla \! _x \cdot \Big((\partial_y u )\nabla \! _x v+ (\partial_y v) \nabla \! _x u\Big) + \partial_y \Big( (\partial_y u)( \partial_y v) - \nabla \! _x u \cdot \nabla \! _x v\Big)  = 0 \label{global_lem_eqn}
\end{equation}
for each $(x,y) \in \Omega$. This can be verified by expanding out the left hand side of \eqref{global_lem_eqn} to find
\[ (\partial_y v)\Delta u+ (\partial_y u) \Delta v, \]
which vanishes in $\Omega$ since both $u$ and $v$ are harmonic in $\Omega$.

\subsection{Irrotational Case}
In this section we proceed formally: we assume $\phi(x,y,t)$ and $\eta(x,t)$ have sufficient decay as $|x|\rightarrow\infty$ for each $(y,t)$ in order for the integrals that follow to exist. Then later, we classify the relevant function spaces so that the integral equations arising in the non-local formulation are well-defined.

Suppose $k \in \mathbf{R}^{n-1}$ and define $v \in C^\infty (\Omega)$ by:
\[ v (x,y) = \exp (-\mathrm{i} k\cdot x + \kappa  y ), \qquad \kappa  = \pm \sqrt{k_1^2 + \cdots+ k_{n-1}^2}.\]
Then $v$ is harmonic in $\Omega$ and the following holds for $(x,y)\in \Omega$:
\[ \nabla \! _x \cdot \left\{ v (-\mathrm{i}k \partial_y \phi + \kappa  \nabla \! _x\phi)\right\} + \partial_y\left\{ v(\kappa  \partial_y \phi + \mathrm{i}k\cdot \nabla \! _x \phi)\right\}  = 0\]
Integrating this expression over $\Omega$ and applying the divergence theorem gives:
\begin{multline} \int_{\mathcal{S}_\eta} e^{-\mathrm{i}k\cdot x +\kappa  y}\left(-\mathrm{i}k\partial_y\phi + \kappa  \nabla \! _x\phi, \kappa  \partial_y\phi + \mathrm{i}k\cdot \nabla \! _x \phi\right) \cdot N_\mathcal{S}\, \dd x \\
 + \int_{\mathcal{B}_h} e^{-\mathrm{i}k\cdot x +\kappa  y}\left(-\mathrm{i}k\partial_y\phi + \kappa  \nabla \! _x\phi, \kappa  \partial_y\phi + \mathrm{i}k\cdot \nabla \! _x \phi\right) \cdot N_\mathcal{B}\, \dd x =0. \label{divthm}
\end{multline}
We have discarded the contribtions from $\partial\mathbf{R}^{n-1}$ using our assumption about the decay of the fields. The contribution from $\mathcal{S}_\eta$ in \eqref{divthm}, is given by the following integral:
\[ \int_{\mathcal{S}_\eta} e^{-\mathrm{i}k\cdot x + \kappa  y}\big[\kappa  \left( \partial_y \phi - \nabla \! _x \phi \cdot \nabla \! _x \eta\right) + \mathrm{i}k\cdot \left( \nabla \! _x \phi + (\partial_y\phi)\nabla \! _x \eta\right)\big]\, \dd x. \]
From \eqref{e3} the first term in the above integrand becomes $\eta_t$, and by \eqref{chain}, the second term becomes $\nabla \! _x q$. This gives the expression
\begin{equation} \int e^{-\mathrm{i}k\cdot x + \kappa  \eta} \left( \kappa  \eta_t + \mathrm{i}k\cdot \nabla \! _x q\right)\, \dd x. \label{contS} \end{equation}
Now we consider the contribution from $\mathcal{B}_h$ in \eqref{divthm}. We introduce $Q(x,t)=\phi (x,-h_0-h,t)$, which represents the potential on the bottom $\mathcal{B}_h$. In analogy with \eqref{chain}, an application of the chain rule yields:
\begin{equation} \nabla \! _x Q = \nabla \! _x \phi -(\partial_y\phi) \nabla \! _x h.\label{chain2}\end{equation}
Recalling that $N_\mathcal{B}=(-\nabla \! _x h, -1)$, the contribution from $\mathcal{B}_h$ in \eqref{divthm} gives
\begin{equation} -\int_{\mathcal{B}_h} e^{-\mathrm{i}k\cdot x + \kappa  y}\big [\kappa  \left( \partial_y \phi +\nabla \! _x \phi \cdot \nabla \! _x h\right) + \mathrm{i}k\cdot \left( \nabla \! _x \phi - (\partial_y\phi)\nabla \! _x h\right)\big ]\, \dd x. \label{bottomterm}\end{equation}
Equation\eqref{e2} implies that $ \nabla\phi \cdot N_\mathcal{B}=0$ on $\mathcal{B}_h$, so the first term in the above integrand vanishes. Using the result in \eqref{chain2} we find that the expression in \eqref{bottomterm} becomes:
\begin{equation} \int e^{-\mathrm{i}k\cdot x - \kappa  (h_0+h)}(-\mathrm{i}k\cdot \nabla \! _x Q)\, \dd x. \label{contB} \end{equation}
Now combining the results in \eqref{contS}, \eqref{contB} and \eqref{divthm} we find that equation \eqref{divthm} reduces to the following equation:
\begin{equation} \int e^{-\mathrm{i}k\cdot x + \kappa  \eta} \left( \kappa  \eta_t + \mathrm{i}k\cdot \nabla \! _x q\right)\, \dd x - \int e^{-\mathrm{i}k\cdot x - \kappa  (h_0+h)}(\mathrm{i}k\cdot \nabla \! _x Q)\, \dd x =0. \label{global_}
\end{equation}
By evaluating \eqref{global_} at $\pm |\kappa |$ and adding/subtracting the resulting equations, we arrive at the following result:
\begin{proposition}[Irrotational Waves]\label{global_inf}
The boundary value problem given in \eqref{e1}-\eqref{e3} is equivalent to the pair of integro-differential equations:
\begin{equation}
\int e^{-\mathrm{i}k\cdot x} \big [ \kappa \eta_t \sinh (\kappa  \eta) + \mathrm{i}k\cdot\nabla \! _x q \cosh (\kappa \eta) - \mathrm{i}k\cdot \nabla \! _x Q \cosh(\kappa (h+h_0))\big] \dd x  = 0 \label{int_temp1}
\end{equation}
and
\begin{equation}
\int e^{-\mathrm{i}k\cdot x} \big[ \kappa  \eta_t \cosh (\kappa  \eta ) + \mathrm{i}k\cdot\nabla \! _x q \sinh (\kappa \eta) + \mathrm{i}k\cdot \nabla \! _x Q \sinh(\kappa (h+h_0))\big] \dd x  = 0 \label{int_temp2}
\end{equation}
valid for each $k\in \mathbf{R}^{n-1}$. These equations, together with the Bernoulli condition \eqref{e5} constitute three equations for the three unknowns $(\eta,q,Q)$.
\end{proposition}
At this point we discuss some rigorous aspects of our approach. We consider appropriate function spaces so that equations \eqref{int_temp1} and \eqref{int_temp2} are well-defined. Of course, one could impose some artificial restrictions on the functions $(\eta,q,Q)$ so that both integrals make sense, but it is important to derive the appropriate function spaces from a more global point of view. For instance, one could impose that both $\nabla \! _x q$, $\nabla \! _x Q$ and $\eta_t$ should be in $L^1(\mathbf{R}^{n-1})$, and that both $\eta$ and $h$ be in $L^\infty(\mathbf{R}^{n-1})$, but there would be no motivation for these constraints in relation to the boundary value problem \eqref{irrotational-n}. Instead, we provide some physically motivated constraints to the problem so that the the integral equations of the non-local approach make sense. We will use some standard results from the theory of Sobolev spaces, and refer the unfamiliar reader to the comprehensive account given in \cite{adams1975sobolev}.

Our starting point is one based on physical grounds: we want to consider solutions to the water wave problem for which the configuration of the system has finite energy. That is to say, we are only interested in water waves for which the total energy:
\[ \iint_\Omega \tfrac{1}{2} |\nabla \phi |^2\, \dd x\, \dd y +\int \left( \tfrac{1}{2} g \eta^2 + \frac{\sigma}{\rho} \left(\sqrt{1+ |\nabla \! _x\eta|^2} -1\right)\right)\dd x, \]
is finite. Clearly this requires $\nabla \phi \in L^2(\Omega)$, and noting
\[ \sqrt{1+t^2} -1 \leq t^2, \qquad t\in \mathbf{R}, \]
we see that it is sufficient that $\eta \in H^1(\mathbf{R}^{n-1})$. With this in mind, we seek to find functions spaces such that the finite energy conditions are satisfied \emph{and} our integral equations are well-defined, given the a priori regularity conditions on $\eta$ and $\nabla \phi$. 

For concreteness, we first set the regularity of the free surface amplitude, $\eta$. The finite energy condition means we need $\eta \in H^s(\mathbf{R}^{n-1})$ for $s\geq 1$. It is inevitable that the regularity of $\eta$ will change the required regularity of the other functions, so we proceed by adjusting $s$ accordingly so that the other estimates can be assessed.

With regards the potential $\phi$: it would be beneficial, from our standpoint, to choose $\nabla\phi$ in a specific function space so that $\nabla \! _x q$ and $\nabla \! _x Q$ are both in $L^1(\mathbf{R}^{n-1})$, owing to the way in which the functions appear in \eqref{int_temp1} and \eqref{int_temp2}. That is to say, we want to classify $X$ so that $\nabla \phi \in X$ guarantees that the \emph{trace} of $\nabla \phi$, i.e. $\nabla \phi|_{\partial\Omega}$, belongs to $L^1(\partial\Omega)$. To ensure the finite energy condition, we would need $X\subset L^2(\Omega)$. Now it is known (e.g. \cite{adams1975sobolev}) that \emph{if} $\partial\Omega$ is sufficiently nice, the trace operator $\mathrm{tr}\, f = f|_{\partial\Omega}$ provides the following continuous map:
\[ \mathrm{tr}\,: H^{r}(\Omega) \rightarrow H^{r-1/2}(\partial\Omega) \qquad \textrm{for $r>\tfrac{1}{2}$}. \]
In this case, ``sufficiently nice'' would be for $\partial\Omega$ to be $C^k$ with $r\leq k$. To enforce this smoothness on the boundary, and hence $\eta$, we must choose $s>\tfrac{1}{2}(n-1)+k$. It would be enough for us to choose $\nabla \phi \in H^{r}(\Omega)$ with $r$ sufficiently large, so that $H^{r-1/2}(\partial\Omega) \subset L^1(\partial\Omega)$. Now the standard embedding result
\begin{equation} H^{m/2}(\mathbf{R}^m) \hookrightarrow L^1(\mathbf{R}^m) \label{embedding}\end{equation}
suggests that we should choose $r\geq\tfrac{n}{2}$. This is enough to make sure our non-local formulation holds.

\begin{lemma}\label{rigorous_estimates_lem}
The global relations \eqref{int_temp1} and \eqref{int_temp2} are well-defined if $\nabla \phi \in H^{s_1}(\Omega)$ and $\eta \in H^{s_1}(\mathbf{R}^{n-1})$, where $s_1 \geq \tfrac{n}{2}$ and $s_2 > n-\tfrac{1}{2}$. In particular, we have the estimates:
\begin{align*} \|\nabla \! _x q\|_{L^1} &\leq c_1  \|\nabla \phi \|_{H^{s_1}} \\ \|\nabla \! _x Q\|_{L^1} &\leq c_2  \|\nabla \phi \|_{H^{s_1}} \\ \|\eta_t \|_{L^1} &\leq c_3  \|\nabla \phi \|_{H^{s_1}},  \end{align*}
where $c_i$, $i=1,2,3$ are constants depending on the dimension $n$ and $\|\eta \|_{H^{s_2}}$.
\end{lemma}
\begin{proof}
Note first that since $s_2>n-\tfrac{1}{2} = \tfrac{1}{2}(n-1) + \tfrac{n}{2}$, it follows that $\eta \in C^{[n/2]}(\mathbf{R}^{n-1})$, where the square bracket denote the integer part. In particular, this tells us that both $\nabla \! _x\eta$ and $\eta$ are bounded, from which we deduce that the global relations \eqref{int_temp1} and \eqref{int_temp2} are well-defined if all of $\eta_t, \nabla \! _xq, \nabla \! _xQ$ are to be in $L^1(\mathbf{R}^{n-1})$. We now establish the estimates in the lemma.

Given $\nabla \phi \in H^{s_1}(\Omega)$, we know that $\mathrm{tr}\, \nabla \phi \in H^{s_1-1/2}(\partial\Omega)$ since our assumption on the regularity of $\eta$, hence $\partial\Omega$, are enough to ensure the trace operator maps $H^{s_1}(\Omega)$ to $H^{s_1-1/2}(\partial\Omega)$. This gives us the following estimates:
\begin{align*} \|\nabla \! _x q\|_{H^{s_1-1/2}} &\leq \|\mathrm{tr}\, \nabla \! _x \phi\|_{H^{s_1-1/2}} + \|\nabla \! _x \eta\, \mathrm{tr}\, \phi_y\|_{H^{s_1-1/2}} \\
 &\leq \|\mathrm{tr}\, \nabla \! _x \phi\|_{H^{s_1-1/2}} + \|\nabla \! _x \eta\|_{L^\infty} \| \mathrm{tr}\, \phi_y\|_{H^{s_1-1/2}} \\
 &\leq \max \{\|\nabla \! _x\eta\|_{L^\infty},1\} \|\mathrm{tr}\, \nabla\phi \|_{H^{s_1-1/2}} \\
 &\leq c_1' \|\nabla\phi \|_{H^{s_1}}.
\end{align*}
The result stated in the lemma follows from the embedding \eqref{embedding} and $s_1 \geq \tfrac{n}{2}$. 

Finally, we see that the regularity of $\eta_t$ is determined through \eqref{irrotational-n}, in particular we have $\eta_t = \mathrm{tr}\, \nabla \phi \cdot N_\mathcal{S}$, i.e.
\[ \eta_t = \mathrm{tr}\, \phi_y - \nabla \! _x\eta \cdot \mathrm{tr} \, \nabla \! _x \phi. \]
Using similar estimates as those for $\nabla \! _x q$, we find the required result.
\end{proof}
Similar arguments can be used for the case of rotational water waves in 2 dimensions (see \S \ref{rotational_section}).

\subsection{The Rotational Two-Dimensional Case}\label{rotational_section}
Let $n=2$ and also confine attention to the case in which $\mathcal{B}_h=\mathcal{B}_0$ is constant. Using a similar approach to the irrotational case, we find the non-local integro-differential equation for the pseudo-potential $\varphi(x,t)$ and the wave height $\eta(x,t)$:
\[ \int_{\mathcal{S}_\eta} e^{-\mathrm{i}kx\mp ky} \left[ \mathrm{i}( \varphi_x + \varphi_y \eta_x) \mp (\eta_x \varphi_x -\varphi_y)\right] \dd x + \int_{\mathcal{B}_0} e^{-\mathrm{i}kx\mp ky}\left[ \mathrm{i}\varphi_x \pm \varphi_y\right]\dd x = 0, \]
which is valid for $k\in \mathbf{R}$. Invoking the boundary conditions \eqref{f2}, \eqref{f3} and \eqref{mess} this expression becomes
\[ \int e^{-\mathrm{i}kx\mp k\eta} \big [ \mathrm{i}\xi_x \mp (\gamma\eta\eta_x - \eta_t)\big ]\dd x + \int e^{-\mathrm{i}kx\pm kh_0} \mathrm{i}\varphi_x (x,-h_0,t)\, \dd x = 0. \]
Subtracting the above two expressions eliminates the $\mathcal{B}_0$ integral completely and we are left with the global relation for two dimensional water waves with constant vorticity $\gamma$:
\begin{equation*} \int e^{-\mathrm{i}kx} \big[ \xi_x \sinh \left( k(\eta+h)\right) - \mathrm{i}(\eta_t - \gamma \eta\eta_x)\cosh (k(\eta+h))\big]\dd x =0, \end{equation*}
where we have dropped the subscript on $h$. This leads us to the following proposition.
\begin{proposition}[Two Dimensional Water Waves with Constant Vorticity]\label{global_inf2}
The boundary value problem in \eqref{f1}-\eqref{f4} for $(\varphi, \eta)$ is equivalent to the following pair of integro-differential equations for $(\xi, \eta)$:
\begin{subequations}\label{g}
\begin{align}
\int e^{-\mathrm{i}kx} \big[ \xi_x \sinh \left( k(\eta+h)\right) - \mathrm{i}(\eta_t - \gamma \eta\eta_x)\cosh (k(\eta+h))\big]\dd x &=0, \label{g1} \\
\xi_t + \tfrac{1}{2} \xi_x^2 -\gamma \partial_x^{-1}\eta_t+ g\eta + \frac{\gamma\eta (2\eta_t \eta_x - 2\xi_x+\gamma\eta)}{2(1+\eta_x^2)} - \frac{(\eta_t + \eta_x \xi_x)^2}{2(1+\eta_x^2)}  &= f(\eta), \label{g2}
\end{align}
\end{subequations}
where $k\in\mathbf{R}$, $\xi (x,t)=\varphi(x,\eta(x,t),t)$ is the pseudo-potential evaluated on $\mathcal{S}_\eta$ and $f(\eta)$ is defined in \eqref{defs}.
\end{proposition}
A similar formulation can be developed if $\mathcal{S}$ is allowed to become multivalued. Indeed, if we assume $\mathcal{S}$ is a differentiable manifold embedded in $\mathbf{R}^2$ via $x=X(\lambda,t)$ and $y=Y(\lambda,t)$, with $|Y|\rightarrow 0$ at $\infty$, then a calculation similar to the previous result gives:
\begin{proposition}[Two Dimensional Water Waves with Constant Vorticity and Multivalued Free Surface]\label{global_inf3}
Let the free surface $\mathcal{S}$ be a 1-dimensional $C^2$-manifold, and let its embedding in $\mathbf{R}^2$ be parameterised by $X(\lambda,t)$ and $Y(\lambda,t)$. Then the boundary value problem in \eqref{f1}-\eqref{f4} for $(\varphi, \mathcal{S})$ is equivalent to the following pair of integro-differential equations for $(\xi,X,Y)$:
\begin{subequations}\label{multig}
\begin{align}
\int e^{-\mathrm{i}kX} \big[ \dot{\xi} \sinh \left( k(Y+h)\right) - \mathrm{i}(\dot{X}Y_t - \dot{Y}X_t - \gamma Y\dot{Y})\cosh (k(Y+h))\big]\dd \lambda &=0, \label{mutlig1} \\
\xi_t + gY -\tfrac{1}{2}X_t^2 - \tfrac{1}{2}Y_t^2 - \gamma \partial_\lambda^{-1} \left( \dot{X} Y_t - \dot{Y}X_t \right)\hspace{6cm} & \nonumber\\ 
+ \frac{\gamma Y(2\dot{X}\dot{Y}Y_t -2\dot{X}^2X_t-2\dot{\xi}\dot{X}+\gamma \dot{X}^2Y)}{2(\dot{X}^2+\dot{Y}^2)} + \frac{(\dot{\xi}-\dot{X}X_t-\dot{Y}Y_t)^2}{2(\dot{X}^2+\dot{Y}^2)} - \sigma \frac{\dot{X}\ddot{Y}-\dot{Y}\ddot{X}}{(\dot{X}^2+\dot{Y}^2)^{3/2}} &=0, \label{multig2}
\end{align}
\end{subequations}
where $k\in\mathbf{R}$, $\xi (\lambda,t)=\varphi(X(\lambda,t),Y(\lambda,t),t)$ is the pseudo-potential evaluated on $\mathcal{S}$.
\end{proposition}

\subsection{Travelling Water Waves}
Let us consider again travelling wave solutions to \eqref{g}, so that $\eta=\eta(x-ct)$ and $\xi=\xi(x-ct)$. Setting $z=x-ct$ and denoting differentiation with respect to $z$ with a prime, the equations \eqref{g} become:
\begin{subequations}\label{g'}
\begin{align}
\int e^{-\mathrm{i}kx} \big[ \xi' \sinh \left( k(\eta+h)\right) + \mathrm{i}(c \eta' + \gamma \eta\eta')\cosh (k(\eta+h))\big]\dd x &=0, \label{g1'} \\
-c \xi' + \tfrac{1}{2} (\xi')^2 + (g+\gamma c)\eta + \frac{\gamma\eta (-2c(\eta')^2 - 2\xi'+\gamma\eta)}{2(1+(\eta')^2)} - \frac{(-c\eta' + \eta' \xi')^2}{2(1+(\eta')^2)}  &= f(\eta). \label{g2'}
\end{align}
\end{subequations}
We recall that that \eqref{g2'} is a quadratic in $\xi'(z)$ which can be solved to give:
\begin{equation} \xi'(z) = \gamma \eta + c - \sqrt{ \left(1+(\eta')^2\right)\left( c^2 + f(\eta)- 2g\eta \right)}, \label{solve_quad}\end{equation}
Replacing in equation \eqref{g1'} $\xi$ by the right hand side of \eqref{solve_quad}, equation \eqref{g1'} becomes a single equation for the unknown $\eta$. In the case of zero vorticity and zero surface tension, the expression simplifies considerably and the global relation \eqref{g1'} becomes:
\[ \int e^{-\mathrm{i} k z} \left[ \Big( 1- \sqrt{\left( 1+(\eta')^2\right) \left( 1- \tfrac{2g}{c^2}\eta \right)}\Big) \sinh (k(\eta+h)) + \mathrm{i}\eta' \cosh (k(\eta+h))\right]\dd z =0. \]

\section{Formal Asymptotic Results for the Two-Dimensional Rotational Case}
In this section we non-dimensionalise the equations in \eqref{g} and approach the problem perturbatively. Throughout this section we make the assumption that each of $\{\eta, \eta_t, \eta_x,\xi_x\}$ are bounded and have sufficient decay at $\infty$ so that the results that follow remain valid. The rigorous justification of these results should be achieved using similar arguments to those in \S 3, but this is not persued here. We suppose $\ell$ is a typical length scale for the wavelengths and $a$ is a typical amplitude of oscillation. Then we make the following substitutions:
\[ x\mapsto \ell x, \quad k\mapsto \frac{k}{\ell}, \quad t\mapsto \frac{\ell }{\sqrt{gh}}t, \quad \xi \mapsto \frac{g\ell a}{\sqrt{gh}} \xi, \quad \eta \mapsto a\eta, \quad \gamma \mapsto \sqrt{\frac{g}{h}} \gamma. \]
We introduce the dimensionless parameters $(\epsilon, \delta)$ defined by:
\[ \epsilon = \frac{a}{h}, \qquad \delta = \frac{h}{\ell}, \]
which are assumed to be small. In this case \eqref{g1} becomes:
\begin{equation}
\int e^{-\mathrm{i}kx} \big\{ \delta^{-1}\xi_x \sinh [\delta k(\epsilon\eta + 1)] - \mathrm{i} [\eta_t - \epsilon  \gamma \eta \eta_x]\cosh [\delta k(\epsilon \eta+1)]\big\}\dd x = 0. \label{int_delta}
\end{equation}
It is straightforward to individually dominate the terms appearing in the integrand, assuming appropriate bounds on $\|\xi_x\|_{L^1 }$, $\|\eta\|_{L^{\infty}}$ and $\|\eta_t\|_{L^1 }$. An application of the dominated convergence theorem allows us to expand the relevant expressions as a power series in $(\epsilon, \delta)$, so \eqref{int_delta} yields the following expression:
\[ \sum_{n,m=0}^\infty \epsilon^n \delta^m \int e^{-\mathrm{i}kx}A_{nm}(k,\eta,\eta_t,\xi)\, \dd x = 0. \]
Using the correspondence between $\partial_x\mapsto \mathrm{i}k$ in the Fourier integral for a \emph{finite} number of terms (so the relevant expression is well-defined in a classical sense) the same equation yields the following:
\[ \sum_{n,m=0}^{\mathrm{finite}} \epsilon^n \delta^m \int e^{\mathrm{i}kx}A_{nm}(\mathrm{i}\partial,\eta,\eta_t,\xi)\, \dd x + \sum_{n,m}^\infty \epsilon^n \delta^m \int e^{\mathrm{i}kx}A_{nm}(k,\eta,\eta_t,\xi)\, \dd x = 0. \]
It is straightforward to bound the terms in the latter integral so the sum is $O(\epsilon^M \delta^N)$ for some specified $M,N$. Using the completeness of the Fourier transform we deduce:
\[ \sum_{n,m=0}^{\mathrm{finite}} \epsilon^n \delta^m A_{nm}(\eta,\eta_t,\xi)\sim 0, \]
where for convenience of notation we have droppped the $\mathrm{i}\partial$ dependence. The first few $A_{nm}$ can be easily computed:
\begin{equation}
A(\eta,\eta_t,\xi) = \left[\begin{array}{llll} \eta_t + \xi_{xx} & 0 & -\tfrac{1}{2}\eta_{txx} - \tfrac{1}{6} \xi_{xxxx} & \cdots \\
                                 (\eta \xi_x)_x- \gamma\eta\eta_x & 0 & -(\eta\eta_t)_{xx} - \tfrac{1}{2}(\eta\xi_x)_{xxx} + \tfrac{1}{2} \gamma (\eta\eta_x)_{xx} & \cdots \\
				0 & 0 & -\tfrac{1}{2} (\eta^2\eta_t)_{xx} - \tfrac{1}{2} (\eta^2\xi_x)_{xxx}+\gamma(\eta^2\eta_x)_{xx} & \cdots \\
				\vdots & \vdots & \qquad\qquad\vdots & \ddots
                               \end{array}\right]. \label{BIG1}
\end{equation}
We note that each of the coefficients is real, and hence \eqref{int_delta} is \emph{one} equation. This is expected, otherwise equations \eqref{g} would constitute an over-determined system of equations for the unknowns $(\eta,\xi)$.

Now we look at the equation \eqref{g2} for $(\eta,\xi)$. Using the non-dimensional parameters this equation becomes:
\begin{multline}
\xi_t + \tfrac{1}{2} \epsilon \xi_x^2  -\gamma \partial_x^{-1} \eta_t + \eta   \\
+ \frac{ \epsilon \gamma \eta(2\epsilon \delta^2 \eta\eta_x-2\xi_x+\gamma \eta)}{2(1+\epsilon^2\delta^2 \eta_x^2)}- \frac{\epsilon \delta^2(\eta_t +\epsilon \eta_x\xi_x)^2}{2(1+\epsilon^2\delta^2\eta_x^2)} - \delta^2 \hat{\sigma} \left( \frac{\eta_x}{\sqrt{1+\epsilon^2\delta^2\eta_x^2}}\right)_x  =0,
\end{multline}
where $\hat{\sigma} \equiv \sigma/gh^2$ is the reciprocal of the bond number. Again we expand the expression in terms the dimensionless parameters $(\epsilon, \delta)$ to find a series of the form:
\[ \sum_{n,m=0}^\mathrm{finite} \epsilon^n \delta^m B_{nm}(\eta, \eta_t, \xi,\xi_t) \sim 0. \]
Again, the computation of the coefficients $B_{nm}$ straightforward:
\begin{equation}
B(\eta,\eta_t,\xi,\xi_t) = \left[\begin{array}{llll} \xi_t+\eta -\gamma \partial_x^{-1} \eta_t & 0 & -\hat{\sigma}\eta_{xx} & \cdots \\
                                 \tfrac{1}{2}\xi_x^2-\gamma \eta \xi_x +\tfrac{1}{2}\gamma^2\eta^2 & 0 &  - \tfrac{1}{2}\eta_t^2 & \cdots \\
				0 & 0 & -\eta_t\eta_x\xi_x & \cdots \\
				\vdots & \vdots & \qquad\vdots & \ddots
                               \end{array}\right]. \label{BIG2}
\end{equation}

The above analysis yields a graded system of equations according to the perturbation parameters $(\epsilon,\delta)$.

\subsection{Solitons and Perturbed Hamiltonian Structure}
In this section we consider the governing equations at each order, using the expressions obtained in \eqref{BIG1} and \eqref{BIG2}. We first investigate the Hamiltonian structure of these equations at each $O(\epsilon^N\delta^M)$, and then demonstrate the existence of solitons at an appropriate truncation. It will be advantageous to introduce the function

\vspace{5mm}
\noindent
\textbf{Order} $\left(\epsilon^0 \delta^0\right)$: To lowest order the evolution equations are:
\begin{subequations}\label{e^0d^0}
\begin{align}
\eta_t &= -\xi_{xx}, \\
\xi_t -\gamma \partial_x^{-1}\eta_t &= -\eta.
\end{align}
\end{subequations}
The right hand side of \eqref{e^0d^0} are in Hamiltonian with respect to $(\eta,\xi)$ with the Hamiltonian density:
\begin{equation*} \mathcal{H}_{00}(\eta,\xi)=\tfrac{1}{2}\int(\eta^2 +\xi_x^2)\,\dd x \end{equation*}
and the standard symplectic structure. That is to say, \eqref{e^0d^0} can be written in the form:
\begin{equation*}
\partial_t \begin{pmatrix} \eta \\ \xi-\gamma\partial_x^{-1}\eta\end{pmatrix} = J\delta \mathcal{H}_{00}, \qquad \textrm{where} \,\,J = \begin{pmatrix} 0 & 1 \\ -1 & 0 \end{pmatrix}\,\, \textrm{and}\,\, \delta = \begin{pmatrix} \delta_\eta \\ \delta_\xi\end{pmatrix}.
\end{equation*}
Here $\delta_\eta$ represents the usual varitational derivative with respect to $\eta$.

\noindent
\textbf{Order} $\left(\epsilon^1 \delta^0\right)$: To next order the equations are:
\begin{subequations}\label{e^1}
\begin{align}
\eta_t &= -\xi_{xx}-\epsilon (\eta \xi_x)_x+ \epsilon\gamma\eta\eta_x, \\
\xi_t-\gamma\partial_x^{-1}\eta_t &= -\eta - \tfrac{1}{2}\epsilon \xi_x^2 + \epsilon \gamma \eta \xi_x-\tfrac{1}{2} \epsilon \gamma^2 \eta^2 .
\end{align}
\end{subequations}
In the same sense, the right hand side of this system is again Hamiltonian with the following Hamiltonian density:
\begin{equation*}\mathcal{H}_{10}(\eta,\xi) = \mathcal{H}_{00}(\eta,\xi) + \tfrac{1}{2}\epsilon\int \left(\eta\xi_x^2+ 2\gamma \eta \eta_x \xi + \tfrac{1}{3}  \gamma^2 \eta^3\right)\,\dd x .\end{equation*}
Indeed \eqref{e^1} can be written in the form:
\begin{equation*}\partial_t \begin{pmatrix} \eta \\ \xi-\gamma\partial_x^{-1}\eta\end{pmatrix} = J\delta \mathcal{H}_{10}.\end{equation*}

\noindent
\textbf{Order} $\left(\epsilon^1 \delta^1\right)$: There is no contribution at this order from the perturbation expansion, so the equations remain unchanged:
\begin{equation*} \partial_t \begin{pmatrix} \eta \\ \xi-\gamma\partial_x^{-1}\eta\end{pmatrix} = J\delta \mathcal{H}_{11}, \end{equation*}
where we have set $\mathcal{H}_{11} = \mathcal{H}_{10} + 0$.

\noindent
\textbf{Order} $\left(\epsilon^0 \delta^2\right)$: At this order these exists a slight complication, because the RHS of the evolution equations involves $\eta_t$:
\begin{equation*} \partial_t \begin{pmatrix} \eta \\ \xi-\gamma\partial_x^{-1}\eta\end{pmatrix} = J\delta \mathcal{H}_{11}  + \delta^2\begin{pmatrix} \tfrac{1}{2} \eta_{txx} + \tfrac{1}{6}\xi_{xxxx} \\ \hat{\sigma}\eta_{xx}\end{pmatrix}.\end{equation*}
However, by using the expression for $\eta_t$ recursively, we can express the RHS in terms of $(\eta,\xi)$ and $x$-derivatives thereof. One must keep track of the order of the relevant terms in the recursive routine. Implementing this approach gives to $O(\delta^2)$ the following equations:
\begin{equation*} \partial_t \begin{pmatrix} \eta \\ \xi-\gamma\partial_x^{-1}\eta\end{pmatrix} = J\delta \mathcal{H}_{11}  + \delta^2\begin{pmatrix}-\tfrac{1}{3}\xi_{xxxx} \\ \hat{\sigma}\eta_{xx}\end{pmatrix}.\label{e^1d^2}\end{equation*}
The right hand side of this system is again Hamiltonian with respect to the Hamiltonian density given by
\begin{equation*}\mathcal{H}_{02}=\mathcal{H}_{11}+\tfrac{1}{2}\delta^2\int (\hat{\sigma}\eta_x^2 - \tfrac{1}{3}\xi_{xx}^2)\, \dd x.\end{equation*}

\noindent
\textbf{Order} $\left(\epsilon^1 \delta^2\right)$: Computing the relevant $O(\epsilon\delta^2)$ terms through the recursion process, we find the following augmented Hamiltonian system:
\begin{equation*} \partial_t \begin{pmatrix} \eta \\ \xi-\gamma\partial_x^{-1}\eta\end{pmatrix} = J\delta \mathcal{H}_{02}  + \epsilon\delta^2\begin{pmatrix}-(\eta\xi_{xx})_{xx} \\ \tfrac{1}{2}\xi_{xx}^2 \end{pmatrix}.\end{equation*}
The associated Hamiltonian density is:
\begin{equation*}\mathcal{H}_{12}=\mathcal{H}_{02}-\tfrac{1}{2}\epsilon\delta^2 \int  \xi_{xx}^2\eta\, \dd x.\end{equation*} 

The above analysis suggests that at each order the system takes the form 
\[ \partial_t \begin{pmatrix} \eta \\ \xi\end{pmatrix} = J\delta \mathcal{H}_{ij} + \begin{pmatrix} 0 \\ \gamma\partial_x^{-1} \eta_t \end{pmatrix} .\]
Using the first of these equations, $\eta_t = \delta_\xi \mathcal{H}_{ij}$, this can be written
\[ \partial_t \begin{pmatrix} \eta \\ \xi\end{pmatrix} = J(\gamma)\delta \mathcal{H}_{ij}, \qquad J(\gamma) = \begin{pmatrix} 0 & 1 \\ -1 & \gamma \partial_x^{-1} \end{pmatrix} \]
where $\mathcal{H}_{ij}$ is an element of the infinite chain of Hamiltonians
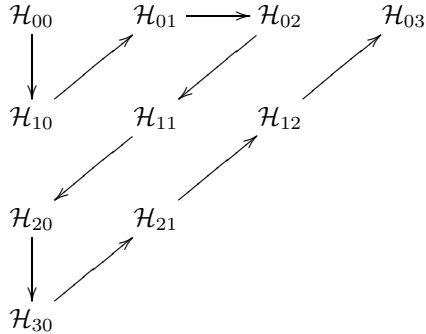
\[ \xymatrix{ 	\mathcal{H}_{00} \ar[d]   & \mathcal{H}_{01} \ar[r]  & \mathcal{H}_{02} \ar[ld] & \mathcal{H}_{03}	\\
		\mathcal{H}_{10} \ar[ru]  & \mathcal{H}_{11} \ar[ld] & \mathcal{H}_{12} \ar[ru] & 	\\
		\mathcal{H}_{20} \ar[d]   & \mathcal{H}_{21} \ar[ru] & 				&	\\
		\mathcal{H}_{30} \ar[ru]  & 			     &				&}	\]
and $J(\gamma)$ represents the perturbed symplectic structure.
\begin{conjecture}
The full system has a perturbed Hamiltonian structure, where the perturbation is defined by the perturbed symplectic form $J(\gamma)$ and the follwing Hamiltonian:
\[ \mathcal{H} = \bigoplus_{n=0}^\infty\bigoplus_{m=0}^\infty \mathcal{H}_{nm} \]
where the grading is with respect to the perturbation parameters $(\epsilon,\delta)$.
\end{conjecture}

It is important to establish the context of the perturbation expansion in this conjecture. Our formulation suggests a perturbed symplectic Hamiltonian structure for the full rotational water wave problem in two dimensions, in the form:
\[ \partial_t \begin{pmatrix} \eta \\ \xi\end{pmatrix} = J(\gamma) \begin{pmatrix} \delta_\eta \\ \delta_\xi\end{pmatrix} \left( \sum_{n,m=0}^\infty  \epsilon^n \delta^m \mathcal{H}_{nm}\right).\]
That is to say, we claim there is a way to ``sum up'' the series above. This is in contrast to the usual perturbative approach to Hamiltonian systems, in which one makes perturbative expansions for the Hamiltonian \emph{and} the symplectic form (e.g. \cite{olver1984hamiltonian}).

We note that our results are in agreement with the work of \cite{constantin2008nearly}, where the authors constructed a nearly Hamiltonian structure for water waves with constant vorticity using the Hilbert transform of the stream function on the free surface. It is well known that the Hilbert transform relates the boundary data for harmonic functions and their conjugates. In light of this, it is not surprising that our formulation is explicit, given that we have chosen to work with the harmonic conjugate of the stream function, rather than the stream function itself.

From this point onwards we truncate the equations at $O(\epsilon^N \delta^M)$, with $N+M>2$. This gives us the perturbed Hamiltonian system:
\begin{align}
\partial_t \begin{pmatrix} \eta \\ \xi \end{pmatrix} &= J(\gamma) \begin{pmatrix} \delta_\eta \\ \delta_\xi \end{pmatrix} \int \tfrac{1}{2}\left[  \left( \eta^2+\xi_x^2\right) +  \epsilon \left(\eta\xi_x^2+ 2\gamma \eta \eta_x \xi + \tfrac{1}{3}  \gamma^2 \eta^3\right) +  \delta^2 (\hat{\sigma}\eta_x^2 - \tfrac{1}{3}\xi_{xx}^2) \right]\, \dd x \nonumber \\
&= \begin{pmatrix} -\xi_{xx}-\epsilon (\eta \xi_x)_x+ \epsilon\gamma\eta\eta_x-\tfrac{1}{3}\delta^2 \xi_{xxxx} \\ -\eta -\gamma \xi_x - \tfrac{1}{2} \epsilon \xi_x^2 + \delta^2 \hat{\sigma} \eta_{xx} - \tfrac{1}{3} \gamma \delta^2 \xi_{xxx}\end{pmatrix} .\label{truncated}
\end{align}
Eliminating $\eta$ from equations \eqref{truncated} and discarding terms of order $O(\epsilon^N \delta^M)$ with $N+M>2$ we find the equation:
\begin{multline} \xi_{tt} - \xi_{xx} +\gamma \xi_{xt} + \epsilon \left( (2+\gamma^2)(\xi_{xt} \xi_x+\gamma \xi_x \xi_{xx}) +(1+\gamma^2)\xi_{xx}\xi_t+ \gamma \xi_t \xi_{xt}\right)\\  + \delta^2 \left( \hat{\sigma}-\tfrac{1}{3}\right) \xi_{xxxx} + \tfrac{1}{3} \gamma \delta^2 \xi_{xxxt} =0. \label{long}
\end{multline}
A formal analysis of the linearisation of \eqref{long} using normal modes shows that solutions will be unstable unless $\hat{\sigma}>1/3$. This is in accordance with the rigorous results for irrotational water waves, e.g. \cite{sun1999net}. To this end we assume $\hat{\sigma}>1/3$ so the corresponding solutions are stable and physically meaningful. We now consider \eqref{long} in two classical scaling regimes and provide a reduction to a nonlinear evolution equation which admits solitons.

\subsection{Soliton Equations}
In this section we consider solitary wave solutions to \eqref{long} in two physically important scaling regimes.

\subsubsection{KdV Scaling Regime with Small Vorticity}
In this case the weak nonlinearity and the weak dispersion balance in the form $\epsilon = \delta^2$, valid for shallow water. We also assume the vorticity is $\mathcal{O}(\epsilon)$, so we make the replacement $\gamma\mapsto \epsilon \gamma$. We use the slow time scale $T=\epsilon t$ and go to a moving frame $X=x-t$. Setting $\xi=\xi(X,T)$ and denoting $\xi_X = \zeta $, the equation \eqref{long} becomes:
\begin{equation}  -2\zeta_T - \gamma \zeta_X - 3 \zeta\zeta_X + \left( \hat{\sigma} - \tfrac{1}{3}\right) \zeta_{XXX} =0. \label{kdvtemp}\end{equation}
This is the celebrated KdV equation. In this case the solitary wave solutions are
\begin{equation} \xi'(z) = (2c-\gamma) \mathrm{sech}^2 \left[ \sqrt{ \frac{\gamma-2c}{\hat{\sigma}-\frac{1}{3}}} \frac{(z-z_0)}{2}\right]. \label{kdvsol}\end{equation}
It is clear that we need $\gamma>2c$ for existence in this case. Setting $\eta = \eta (z)$ in \eqref{truncated} and integrating up we find the following quadratic in $\eta$:
\[ \tfrac{1}{2} \gamma (\epsilon\eta)^2 + (c-\xi') \epsilon \eta - (\xi' + \tfrac{1}{3} \xi ''') =0. \]
In the case $\gamma \neq 0$, we find the following two solutions:
\begin{equation} \eta(z) = \frac{\xi' - c \pm \sqrt{ (\xi'-c)^2 + 2\gamma (\xi' + \tfrac{1}{3} \xi''')}}{\epsilon \gamma} \label{solitons} \end{equation}
of which only the $(+)$ solution obeys $\eta\rightarrow 0$ as $|z|\rightarrow \infty$. These solutions are valid when the term in the square root is positive and the steepness of the corresponding wave increases as this term approaches zero. It is clear from \eqref{solitons} these solutions break down if
\begin{equation} \left( \frac{\dd \xi}{\dd z} -c\right)^2 + 2\gamma \left( \frac{\dd \xi}{\dd z} + \frac{1}{3} \frac{\dd^3 \xi}{\dd z^3}\right) \leq 0 \label{temp} \end{equation}
for some value of $z$. Using \eqref{kdvsol}, this results in a transcendental equation for $z$. However, note that \eqref{kdvtemp}, after one integration, gives $\xi'''$ in terms of $\xi'$. Consequently \eqref{temp} becomes the following quadratic in $\xi'$:
\[ \left[ 1+ \frac{\gamma}{\hat\sigma - \tfrac{1}{3}}\right] (\xi')^2 + 2 \left[ \gamma-c + \frac{\gamma (\gamma-2c)}{3(\hat{\sigma}-\tfrac{1}{3})}\right] \xi' +c^2 \leq 0 \]
By examining the discriminant of this quadratic, one can identify a \emph{subset} of the parameter space $(\hat{\sigma},c,\gamma)$ for which the solution \eqref{solitons} are well-defined.

\subsubsection{Weak Nonlinear-Dispersion Balance}
In the balance between nonlinearity (governed by $\epsilon \ll 1$) and dispersion ($\delta\ll 1$) one assumes the two fundamental length scales are comparable, so we set $\epsilon =\delta$. We then operate with the fast scales $\epsilon X=x$, $\epsilon T=t$. In this regime the equation \eqref{long} transforms to:
\begin{multline} \xi_{TT} - \xi_{XX} +\gamma \xi_{XT} + (2+\gamma^2)(\xi_{XT} \xi_X+\gamma \xi_X \xi_{XX}) +(1+\gamma^2)\xi_{XX}\xi_T+ \gamma \xi_T \xi_{XT}\\  +  \left( \hat{\sigma}-\tfrac{1}{3}\right) \xi_{XXXX} + \tfrac{1}{3} \gamma  \xi_{XXXT} =0 \label{long2}
\end{multline}
For solitons, we look for travelling wave solutions of the form $\xi=\xi(X-cT)$. Using this ansatz in \eqref{long2} we get:
\begin{equation} (c^2-1-\gamma c)\xi'' +   ((2+\gamma^2)(\gamma-c)-c(1+\gamma^2)+c^2\gamma) \xi' \xi'' + \left( \hat{\sigma}-\tfrac{1}{3}-\tfrac{1}{3}\gamma c\right) \xi'''' =0 \label{ode} \end{equation}
where the prime denotes differentiation with respect to $z=X-cT$. Integrating up the equation in \eqref{ode} we find the classical $\mathrm{sech}^2$ soliton solution:
\begin{subequations}\label{xi}
\begin{equation} \xi'(z) = A(\gamma,c)\, \mathrm{sech}^2 \left[ \sqrt{ \frac{ c^2-(1+\gamma c)}{\tfrac{1}{3}(1+\gamma c)-\hat{\sigma}}} \frac{(z-z_0)}{2}\right],  \end{equation}
where the amplitude is given by
\begin{equation}  A(\gamma,c) = 3 \left[ \frac{ 1+\gamma c - c^2}{(2+\gamma^2)(\gamma-c)-c(1+\gamma^2)+c^2\gamma}\right]. \end{equation}
\end{subequations}
The corresponding expression for the amplitude is given by \eqref{solitons}. As with the KdV scaling regime, we can identify a subset of the parameter space for which these solutions are well-defined. The condition in \eqref{temp} is now suplimented with \eqref{ode} which, after one integration, gives $\xi'''$ in terms of $\xi'$ and $(\xi')^2$. Examining the discriminant of the resulting quadratic gives sufficient conditions for the amplitude to be well-defined.

\section{A Rigorous Derivation of the Linear Limit of Irrotational Water Waves}
Here we concentrate on the nonlinear boundary value problem described by \eqref{e1}-\eqref{e3} and \eqref{e5} in the case $\mathcal{B}_h \equiv\mathcal{B}_0$ is constant. We are interested in the asymptotic reduction of these equations to the linear limit, given certain smallness conditions on the norms of the functions in the relevant function spaces.

For ultimate flexibility we work with the space of tempered distributions $\mathcal{S}'(\mathbf{R}^{n-1})$, on which the Fourier transform:
\[ \mathcal{F}_{x\mapsto k}: f \mapsto \hat{f}(k) = \int e^{-\mathrm{i} k\cdot x} f(x)\, \dd x \] 
is defined using the duality with $\mathcal{S}(\mathbf{R}^{n-1})$, the space of Schwartz functions: smooth functions with rapid decay.

Evaluating the integro-differential equations in Proposition \ref{global_inf} for $h=0$, i.e $\mathcal{B}_h$ constant, we find:
\begin{equation}
\int e^{-\mathrm{i}k\cdot x} \big[ \kappa \eta_t \sinh (\kappa  \eta) + \mathrm{i}k\cdot\nabla \! _x q \cosh (\kappa \eta) - \mathrm{i}k\cdot \nabla \! _x Q \cosh(\kappa h_0)\big] \dd x  = 0 \label{er1}
\end{equation}
and
\begin{equation}
\int e^{-\mathrm{i}k\cdot x} \big[ \kappa  \eta_t \cosh (\kappa  \eta ) + \mathrm{i}k\cdot\nabla \! _x q \sinh (\kappa \eta) + \mathrm{i}k\cdot \nabla \! _x Q \sinh(\kappa h_0)\big] \dd x = 0. \label{er2}
\end{equation}
Now multiplying \eqref{er1} by $\sinh (\kappa  h_0)$ and \eqref{er2} by $\cosh (\kappa h_0)$ and adding, we find the following equations:
\begin{equation}
\int e^{-\mathrm{i}k\cdot x} \big[ \kappa \eta_t \cosh (\kappa  (\eta+h_0)) + \mathrm{i}k\cdot\nabla \! _x q \sinh (\kappa (\eta+h_0)) \big] \dd x = 0. \label{grlin}
\end{equation}
This integro-differential equation is valid for $k\in\mathbf{R}^{n-1}$ and constitutes the global relation for the problem in \eqref{e1}-\eqref{e3} in the particular case where $\mathcal{B}_h$ is constant. Our aim is to make suitable approximations in \eqref{grlin} and to bound the relevant errors.

The linear limit is found by assuming $(\eta,\nabla \! _xq)$ and certain derivatives thereof are small, in an appropriate sense, and discarding terms that are smaller. We work on $ \mathcal{S}'(\mathbf{R}^{n-1})$ with the following assumptions:
\[ \|\eta\|_{L^\infty}<\epsilon,\qquad \|\eta_t \|_{L^1 } < \epsilon, \qquad  \|\nabla \! _x q\|_{L^1 } < \epsilon, \]
for small $\epsilon$. These estimates can be realised by restricting $\|\nabla\phi\|_{H^{s_1}(\Omega)}$ and $\|\eta \|_{H^{s_2}}$ are sufficiently small, and using arguments analogous to those used in the proof of lemma \ref{rigorous_estimates_lem}.

We now concentrate solely on the first term in the first integral in \eqref{grlin}, since the results for the other two terms can be derived analogously.
\begin{lemma}\label{bound}
Let $\hat{\eta}$ denote the Fourier transform of $\eta \in \mathcal{S}'(\mathbf{R}^{n-1})$ and let the norms $\|\nabla \phi\|_{H^{s_1}}$ and $\|\eta\|_{H^{s_2}}$ appearing in Lemma \ref{rigorous_estimates_lem} be such that::
\[ \max\{\|\eta _t\|_{L^1},\|\eta\|_{L^\infty}\}< \epsilon.\]
Then the following estimate holds:
\[ \int e^{-\mathrm{i}k\cdot x}  \frac{\eta_t \cosh (\kappa  (\eta+h_0))}{\cosh (\kappa  (\epsilon+h_0))}\, \dd x = \hat{\eta}_t + O(\epsilon^2), \]
valid for $\kappa < O(1/\epsilon)$.
\end{lemma}
\begin{proof}
We breifly outline how to prove the following basic estimate:
\[ \left| \int e^{-\mathrm{i}k\cdot x}  \frac{\eta_t \cosh (\kappa  (\eta+h_0))}{\cosh\big(\kappa (\epsilon+h_0)\big)}\, \dd x - \hat{\eta}_t  \right| < 2\left( 1-e^{-2\kappa \epsilon}\right) \epsilon.  \]
First we note that the LHS can be written as:
\begin{equation} \frac{1}{\cosh\big(\kappa (\epsilon+h_0)\big)} \left| \int e^{-\mathrm{i}k\cdot x} \big[ \eta_t \cosh (\kappa (\eta+h_0)) - \eta_t \cosh (\kappa (\epsilon+h_0))\big] \dd x\right|. \label{bd1}\end{equation}
Given that $\|\eta\|_{L^\infty}<\epsilon$, the following identity holds almost everywhere in $\mathbf{R}^{n-1}$:
\[ \cosh \big(\kappa (\eta+h_0)\big) - \cosh \big(\kappa (\epsilon+h_0)\big) = e^{\kappa (h_0 + \epsilon)} \frac{\big(1-e^{-\kappa|\eta-\epsilon|}\big)}{2} + e^{-\kappa (h_0 + \epsilon)}\frac{ \big(1-e^{-\kappa|\eta+\epsilon|}\big)}{2} \]
Using this in \eqref{bd1} we find:
\begin{align*} \left| \int e^{-\mathrm{i}k\cdot x}  \frac{\eta_t \cosh (\kappa  (\eta+h_0))}{\cosh \big(\kappa (\epsilon+h_0)\big)}\, \dd x - \hat{\eta}_t  \right| &\leq \frac{e^{\kappa(h_0+\epsilon)}}{\cosh \big(\kappa (h_0+\epsilon)\big)} \int |\eta_t|\big( 1-e^{-\kappa |\eta-\epsilon|}\big)\, \dd x \\
&\leq 2 \big(1-e^{-2\kappa \epsilon}\big) \|\eta_t\|_{L^1}.
\end{align*}
Using the a priori bound $\|\eta_t\|_{L^1}<\epsilon$, the result follows.
\end{proof}
\begin{remark}The bound in lemma \ref{bound} can be made sharper if $\eta_t \in H^s(\mathbf{R}^{n-1})$ ($s>1$) with $\|\eta_t\|_{H^s}<\epsilon$, and then use integration by parts $s$ times. However, to do this we would require more restrictive estimates for $\|\nabla\phi\|$ and $\|\eta\|$. This would improve our bound by some algebraic order in $k$, but the bound is not uniform in $\kappa $. The asymptotic estimate in lemma \ref{bound} is valid for $\kappa <O(1/\epsilon)$, which means that the estimate is valid for sufficiently large wavelengths. To understand the shorter wave length dynamics, we would need more knowledge of the smoothness of the functions appearing in \eqref{grlin}. If the functions are smooth, then the rapid oscillations as $|k|\rightarrow \infty$ owing to the $\exp({\mathrm{i} k\cdot x})$ term would counter the growth and global estimates could be established.
\end{remark}

The second term in \eqref{grlin}, i.e the term
\[ \int e^{-\mathrm{i}k\cdot x} k\cdot\nabla \! _x q \sinh (\kappa (\eta+h_0)) \, \dd x, \]
can be estimated in an entirely analagous fashion with the previous result. Indeed, using a similar argument and the a priori bound $\|q_x\|_{L^1}<\epsilon$ we find:
\begin{equation}
\left| \int e^{-\mathrm{i}k\cdot x}  \left(\frac{\mathrm{i}k\cdot \nabla \! _x q \sinh (\kappa  (\eta+h_0))}{\cosh (\kappa  (\epsilon+h_0))}- \mathrm{i}k\cdot \nabla \! _x q \tanh[\kappa  h_0]\right)\, \dd x  \right| < 2\left( 1-e^{-2\epsilon \kappa }\right) \epsilon.
\end{equation}
Then using $\nabla \! _x \mapsto \mathrm{i}k$ in the Fourier integral, we arrive at the following result:
\begin{lemma}\label{bound2}
Let $\hat{q}$ denote the Fourier transform of $q \in \mathcal{S}'(\mathbf{R}^{n-1})$ and let the norms $\|\nabla \phi\|_{H^{s_1}}$ and $\|\eta\|_{H^{s_2}}$ appearing in Lemma \ref{rigorous_estimates_lem} be such that::
\[ \max\{\|\eta \|_{L^\infty},\|\nabla \! _xq\|_{L^1}\}< \epsilon.\]
Then the following estimate holds:
\[ \int e^{-\mathrm{i}k\cdot x} \left[ \frac{\mathrm{i}k\cdot \nabla \! _x q \sinh (\kappa  (\eta+h_0))}{\cosh (\kappa  (h_0+\epsilon))}\right]\, \dd x = - \kappa  \tanh [\kappa  h_0] \hat{q} + O(\epsilon^2), \]
valid for $\kappa < O(1/\epsilon)$.
\end{lemma}
The results in lemmas \ref{bound} and \ref{bound2} provide us with a rigorous linear reduction of \eqref{grlin} valid for sufficiently long wave lengths, as summarised in the following.
\begin{proposition}\label{linear11}
Let $(q,\eta)$ satisfy the boundary value problem in \eqref{e1}-\eqref{e3} for the case of the flat bottom $\mathcal{B}_{h}=\mathcal{B}_0$. Let the norms $\|\nabla \phi\|_{H^{s_1}}$ and $\|\eta\|_{H^{s_2}}$ appearing in Lemma \ref{rigorous_estimates_lem} be such that:
\[ \max \{  \|\eta_t\|_{L^1}, \|\eta\|_{L^\infty}, \|\nabla \! _x q\|_{L^1} \}<\epsilon. \]
Then the following estimate is valid
\begin{equation}  \hat{\eta}_t - \kappa  \tanh[\kappa  h_0] \hat{q}  = O(\epsilon^2), \label{linear1}\end{equation}
for $\kappa < O(1/\epsilon)$.
\end{proposition}
\begin{proof}
As a consequence of the global relation \eqref{grlin}, we observe the identity:
\begin{multline} \hat{\eta}_t - \kappa  \tanh[\kappa  h_0] \hat{q} = \left(\hat{\eta}_t - \int e^{-\mathrm{i}k\cdot x}  \frac{\eta_t \cosh (\kappa  (\eta+h_0))}{\cosh (\kappa  (h_0+\epsilon))}\, \dd x\right) \\ - \left(\kappa  \tanh[\kappa  h_0] \hat{q} +  \int e^{-\mathrm{i}k\cdot x} \left[ \frac{\mathrm{i}k\cdot \nabla \! _x q \sinh (\kappa  (\eta+h_0))}{\cosh (\kappa  (h_0+\epsilon))}\right] \dd x\right),\end{multline}
which holds for all $k\in \mathbf{R}^{n-1}$. The result follows from an application of lemmas \ref{bound} and \ref{bound2}.
\end{proof}
In what follows we will use the same contraints on $(\eta,q)$ and linearise the Bernoulli equation \emph{in the Fourier space}, so that we find an additional equation that couples $\eta$ and $q$, or rather $\hat{\eta}$ and $\hat{q}$.

Recall that the Bernoulli condition reads
\begin{equation} q_t + \tfrac{1}{2}|\nabla \! _x q|^2 + g\eta - \frac{(\eta_t + \nabla \! _x\eta \cdot\nabla \! _x q)^2}{2(1+|\nabla \! _x \eta|^2)} = \frac{\sigma}{\rho} \nabla \! _x \cdot \left( \frac{ \nabla \! _x\eta}{\sqrt{1+|\nabla \! _x\eta|^2}}\right). \label{bernoulli}\end{equation}
It is convenient to rewrite \eqref{bernoulli} as follows
\begin{equation} q_t + g\eta - \frac{\sigma}{\rho} \Delta\eta + N(q,\eta) = 0, \label{split}
\end{equation}
where $\Delta$ is the standard Laplacian on $\mathbf{R}^{n-1}$ and $N(q,\eta)$ is defined as
\begin{equation} N(q,\eta) \defn \tfrac{1}{2}|\nabla \! _xq|^2+\frac{\sigma}{\rho} \nabla \! _x \cdot\left( \frac{ \nabla \! _x\eta}{\sqrt{1+|\nabla \! _x\eta|^2}} - \nabla \! _x \eta\right) - \frac{(\eta_t + \nabla \! _x\eta \cdot\nabla \! _x q)^2}{2(1+|\nabla \! _x \eta|^2)}. \end{equation}
By applying the Fourier transform to \eqref{split} it can be shown that under the assumption that $\|\eta\|_{H^2  }$ and $\|\nabla \! _x q\|_{L^2 }$ are sufficiently small, $\hat{N}(q,\eta)$ is negligable so that the linear terms in \eqref{split} are a good approximation for the dynamics. These estimates on $\|\eta\|$ and $\|\nabla \! _x q\|$ can be obtained by choosing $\|\nabla \phi\|_{H^{s_1}}$ and $\|\eta\|_{H^{s_2}}$ of lemma \ref{rigorous_estimates_lem} are sufficiently small. For example, $\|\nabla \! _x q\|_{L^2}$ is dominanted since $s_1-\tfrac{1}{2}
 \geq \tfrac{1}{2}(n-1)$ so that
\[ H^{s_1-1/2}(\mathbf{R}^{n-1}) \hookrightarrow L^2(\mathbf{R}^{n-1}). \]
Using similar estimates to those in lemmas \ref{bound} and \ref{bound2} (see appendix) we arrive at the following result.
\begin{proposition}\label{linear22}
Let $(q,\eta)$ satisfy the Bernoulli condition \eqref{bernoulli}. Let the norms $\|\nabla \phi\|_{H^{s_1}}$ and $\|\eta\|_{H^{s_2}}$ appearing in Lemma \ref{rigorous_estimates_lem} be such that:
\[ \max \{ \|\eta_t\|_{L^2 }, \|\eta\|_{H^2  }, \|q_x\|_{L^2 }  \} < \epsilon. \]
Then the following estimate is valid
\begin{equation}
\hat{q}_t + g\hat{\eta} + \frac{\sigma}{\rho} \kappa^2 \hat{\eta} = O(\epsilon^2 ) \label{linear2}
\end{equation}
uniformly in $\kappa$.
\end{proposition}
\begin{remark}
The results from Propositions \ref{linear11} and \ref{linear22} yield the classical dispersion relation for linearised water waves. Indeed, differentiating \eqref{linear1} with respect to $t$ and using \eqref{linear2} we find the following equation for $\hat{\eta}$:
\[ \hat{\eta}_{tt} + \kappa g \tanh[\kappa h_0]\left( 1 + \frac{\sigma}{g\rho}\kappa^2\right)\hat{\eta} = 0, \]
where we have discarded the $O(\epsilon^2)$ terms.
\end{remark}

\section{Conclusions}
We have presented the following generalisations of the results of \cite{ablowitz2006nnl}: (a) We have extended and rigorously described the linearisation procedure; (b) we have considered the case of constant vorticity; and (c) we have incorperated the effect of a multi-valued free surface. Also, following \cite{tio2010}, we have derived an upper bound for the wave height in the case of constant vorticity and non-zero surface tension.

We recall that \cite{craig1993nsg} introduced an elegant Dirichlet to Neumann operator $G(\eta)$ associated with the velocity potential and also obtained a series for the operator $G(\eta)$, valid for small $\eta$. The integral equations presented in \cite{ablowitz2006nnl} and here, can be considered as the \emph{summation} of the above series, i.e. the series of \cite{craig1993nsg} is the Neumann series of the integral equations derived by \cite{ablowitz2006nnl} and here.

It appears that the new formulation provides an alternative, perhaps simpler approach, for (a) the numerical investigation of water waves; (b) the derivation of various asymptotic limits; (c) the rigorous analysis of water waves. 

Regarding (a), we recall that two-dimensional lumps were computed by \cite{ablowitz2006nnl} in the case of sufficient surface tension.

Regarding (b) we recall that various asymptotic equations, including the Boussinesq, Benney-Luke and the nonlinear Schr\"odinger equations, were derived by \cite{ablowitz2006nnl}. Similarly, an appropriate Boussinesq type equation in the case of constant vorticity is found in \S 4. We observe that these equations have been derived by several other authors, however, it appears that the new formulation provides a straightforward way of deriving these equations. For example, regarding our results in \S 4, we note that \cite{chow1989sos} uses a perturbative approach which involves solving many PDEs arising from the consistency of the perturbation expansion with the boundary conditions. To solve these PDEs, the author confines attention to seperable solutions. The entire paper is devoted to computing the first few coefficients of the perturbation expansion. The author comments that this method is unweildy at anything beyond second order. Our paper gives a more direct approach, and can easily be extended to higher orders.

Regarding (c) we have shown in \S 5 that standard PDE techniques can be used for the rigorous analysis of water waves, at least in the linear limit. The new formulation suggests a rigorous methodology which differs drastically from that employed in the important works of \cite{amick1982scw, constantin2004esp,wu1997wps}; the extension of the results of \S 5 to the nonlinear problem is a work in progress.

Finally, we note that in the classical works of \cite{benjamin1962sws}, the author does not include the effect of surface tension. In contrast, our formulation also includes the effect of surface tension. Furthermore, our formalism can easily be extended to the three dimensional case. The main advantages of the new approach presented by \cite{ablowitz2006nnl} and in this paper, are a consequence of the \emph{explicit} nature of the equations in Propositions \ref{global_inf}, \ref{global_inf2} and \ref{global_inf3}.

\subsection*{Acknowledgements}
The first author is grateful for the support of Emmanuel College, Cambridge. The second author acknowledges partial support from the Guggenheim Foundation, USA.

\section{Appendix}
Here we prove the result in Proposition \ref{linear22}. It suffices to prove the estimate:
\[ \|N(\eta,q)\|_{L^1} = O(\epsilon^2) \]
given that $\max\{ \|\eta\|_{H^2}, \|q\|_{H^1}\|\}<\epsilon$. Recall that:
\[ N(q,\eta) \defn \tfrac{1}{2}|\nabla \! _xq|^2+\frac{\sigma}{\rho} \nabla \! _x \cdot\left( \frac{ \nabla \! _x\eta}{\sqrt{1+|\nabla \! _x\eta|^2}} - \nabla \! _x \eta\right) - \frac{(\eta_t + \nabla \! _x\eta \cdot\nabla \! _x q)^2}{2(1+|\nabla \! _x \eta|^2)}. \]
The following estimate is clear from the definition:
\[ \|e^{\mathrm{i} k\cdot x}N(\eta,q)\|_{L^1} \leq \tfrac{1}{2} \|\nabla \! _x q\|_{L^2 }^2 + \tfrac{\sigma}{\rho} \left\| \nabla \! _x \cdot \left( \tfrac{\nabla \! _x \eta}{\sqrt{1+|\nabla \! _x\eta|^2}}-\nabla \! _x \eta\right)\right\|_{L^1} + \tfrac{1}{2} \left\| \tfrac{\eta_t + \nabla \! _x \eta\cdot \nabla \! _x q}{\sqrt{1+|\nabla \! _x \eta|^2}}\right\|^2_{L^2 }.\]
The first term is clearly $O(\epsilon^2)$ since $\|q\|_{H^1}<\epsilon$. For the third term, we use the estimate
\begin{align}
\left\| \tfrac{\eta_t + \nabla \! _x \eta\cdot \nabla \! _x q}{\sqrt{1+|\nabla \! _x \eta|^2}}\right\|^2_{L^2 } &\leq \left( \|\eta_t\|_{L^2 } + \left\| |\nabla \! _x q|\left[ \tfrac{ |\nabla \! _x \eta| }{\sqrt{1+|\nabla \! _x\eta|^2}}\right]\right\|_{L^2 } \right)^2 \nonumber \\
&\leq \left( \|\eta_t\|_{L^2 } + \|\nabla \! _xq \|_{L^2 } \right)^2 \label{a2}
\end{align}
which gives the required $O(\epsilon^2)$ bound. For the second term, note that
\begin{multline*} \left\|\nabla \! _x \cdot \left( \frac{\nabla \! _x \eta}{\sqrt{1+|\nabla \! _x\eta|^2}}-\nabla \! _x \eta\right)\right\|_{L^1} \leq \left\|\left( \frac{1}{\sqrt{1+|\nabla \! _x\eta|^2}}-1\right)\Delta \eta\right\|_{L^1}\\+  \sum_{i,j} \left\| \frac{\eta_i \eta_j \eta_{ij}}{(1+|\nabla \! _x\eta|^2)^{3/2}}\right\|_{L^1} \end{multline*}
where $\eta_i \equiv \partial \eta /\partial x_i$. Rewriting the second term, we observe the following estimate:
\begin{align} \sum_{i,j} \left\| \frac{\eta_i \eta_j \eta_{ij}}{(1+|\nabla \! _x\eta|^2)^{3/2}}\right\|_{L^1} &= \sum_{i,j}\left\| \left( \frac{\eta_i}{(1+\sum_k\eta_k^2)^{3/2}}\right)\eta_{ij}\eta_j\right\|_{L^1}\nonumber \\
 &\leq \sum_{i,j}\left\| \left(\frac{\eta_i}{(1+\eta_i^2)^{3/2}}\right)\eta_{ij}\eta_j\right\|_{L^1}\nonumber \\
 &\leq \sum_{i,j}\tfrac{2}{3\sqrt{3}}\left\| \eta_{ij}\eta_j\right\|_{L^1} \nonumber \\
 &\leq \sum_{i,j}\tfrac{2}{3\sqrt{3}}\left\| \eta_{ij}\|_{L^2 }\|\eta_j\right\|_{L^2 }, \label{a0}
\end{align}
where we applied the Cauchy-Schwarz inequality. Similarly, we find
\begin{equation} \left\|\left( \frac{1}{\sqrt{1+|\nabla \! _x\eta|^2}}-1\right)\Delta \eta\right\|_{L^1} \leq \tfrac{1}{2} \|\nabla \! _x\eta\|_{L^2 } \|\Delta \eta\|_{L^2 }\label{a1}. \end{equation}
Combining \eqref{a0} and \eqref{a1} we find
\begin{equation} \left\|\nabla \! _x \cdot \left( \frac{\nabla \! _x \eta}{\sqrt{1+|\nabla \! _x\eta|^2}}-\nabla \! _x \eta\right)\right\|_{L^1} \leq \|\eta\|_{H^2}^2. \label{a3}\end{equation}
From the estimates in \eqref{a2} and \eqref{a3} it is follows that $\|N(\eta,q)\|_{L^1}=O(\epsilon^2)$.

\nocite{ablowitz2009asymptotic}
\nocite{haut2009reformulation}

\bibliographystyle{plain}

\end{document}